\documentclass[a4paper,12pt]{article}
\usepackage[T1]{fontenc}
\usepackage{latexsym}
\usepackage{amssymb,amsmath,mathtools}
\usepackage{graphics}
\usepackage[a4paper]{geometry}
\newtheorem{thm}{Theorem}

\newtheorem{prop}[thm]{Proposition}

\newtheorem{defn}[thm]{Definition}
\newcommand{\C}{\mathbb{C}}

\newcommand{\Z}{\mathbb{Z}}
\newcommand{\N}{\mathbb{N}}

\newcommand{\Tr}{\operatorname{Tr}}
\newcommand{\ep}{\varepsilon}

\title{\textbf{Right submodules of finite rank for von Neumann dynamical systems}}
\author{Paul Jolissaint}

\begin{document}

\maketitle

\begin{abstract}
Let $(M,\tau,\sigma,\Gamma)$ be a (finite) von Neumann dynamical system and let $N$ be a unital von Neumann subalgebra of $M$ which is $\Gamma$-invariant. If $V\subset L^2(M)$ is a right $N$-submodule whose projection $p_V$ has finite trace in $\langle M,e_N\rangle$ and is $\Gamma$-invariant, then we prove that, for every $\varepsilon>0$, one can find a $\Gamma$-invariant submodule $W\subset V$ which has finite rank and such that $\Tr(p_V-p_W)<\varepsilon$. In particular, this answers Question 4.2 of \cite{AET}.
\par\vspace{3mm}\noindent
\end{abstract}

\section{Introduction}

Let $(M,\tau,\sigma,\Gamma)$ be a von Neumann dynamical system, \emph{i.e.} $M$ is a finite von Neumann algebra endowed with
\begin{itemize}
	\item a normal, finite, faithful, normalized trace $\tau$,
	\item an action $\sigma$ of a countable group $\Gamma$ such that every $\sigma_g$ is $\tau$-preserving for every $g\in \Gamma$.
\end{itemize}
In order to prove an extention to von Neumann dynamical systems of multiple recurrence Szemer\'edi's Theorem, the authors of \cite{AET} need the following lemma (\cite{AET}, Lemma 4.1, stated here somewhat vaguely): 
\par
\emph{Let $(M,\tau,\alpha)$ be a finite von Neumann algebra endowed with a $\tau$-preserving action of $\Z: n\mapsto\alpha^n$. If $N$ is an $\alpha$-invariant von Neumann subalgebra of the center of $M$ and if $V\subset L^2(M)$ is an invariant right-$N$-submodule with finite trace, then one can find an invariant $N$-submodule $W\subset V$ of finite rank $r$, arbitrarily close to $V$ in an appropriate sense, and the associated action of $\alpha$ on $W$ is described by a suitable unitary $u\in U(M_r(N))$
(all relevant definitions are recalled below).}
\par
The proof rests on the decomposition of $L^2(M)$ in a direct integral based on $N$, and that is why $N$ is assumed to lie in the center of $M$. Guessing that the conclusions of the lemma should hold under more general hypotheses, the authors of \cite{AET} ask if it is indeed the case in Remark 4.2 of their article. 
\par\vspace{3mm}
Thus, the aim of the present notes is to provide more general statements of the above lemma. As a matter of fact, we present two variants. In the first one (Proposition 1), we assume that $N$ is at the opposite of being central, namely that $N$ is a type II$_1$ factor, and we get a neater version than in the general case which is the subject of Proposition 2. The first result and its proof are contained in Section 2, and the second one is treated in the last section. We recall below the main definitions and results that will be needed in the forthcoming sections. 

\par\vspace{3mm}
Standard notation in the framework of finite von Neumann algebras will be used; they are borrowed from A. Sinclair's and R. Smith's monograph
\cite{SS}. Let $(M,\tau)$ be as above, and let $L^2(M)=L^2(M,\tau)$ be the  associated Hilbert space given by GNS construction. It is an $M$-bimodule, and $M$ embeds into $L^2(M)$ as a dense subspace. We denote by $\xi$ the image of $1\in M$ in $L^2(M)$, so that $\tau(x)=\langle x\xi,\xi\rangle$ for every $x\in M$. 
Let $1\in N\subset M$ be a von Neumann subalgebra of $M$. We denote by $E_N$ the $\tau$-preserving conditional expectation of $M$ onto $N$, and by $e_N$ the extention of $E_N$ to $L^2(M)$;
$\langle M,e_N\rangle=JN'J$ denotes the associated basic construction algebra.
\par\vspace{3mm}
Recall also that $L^1(M)$ is the completion of $M$ with respect to the norm 
$$
\Vert x\Vert_1:=\tau(|x|)=\sup\{|\tau(xa)|: a\in M,\Vert a\Vert\leq 1\}.
$$
It is an $M$-bimodule because it is easy to check that
$\Vert axb\Vert_1\leq \Vert a\Vert\Vert x\Vert_1\Vert b\Vert$ for all $a,x,b\in M$. Furthermore, the mapping $(x,y)\mapsto xy$ defined from $M\times M$ to $M$ extends to a mapping from $L^2(M)\times L^2(M)$ to $L^1(M)$ because $\Vert xy\Vert_1\leq\Vert x\Vert_2\Vert y\Vert_2$ by Cauchy-Schwarz inequality.
In particular, setting $y=1$, we get the inequality $\Vert x\Vert_1\leq\Vert x\Vert_2$ which implies that $L^2(M)\subset L^1(M)$.
Similarly, the conditional expectation $E_N$ extends to a contraction from $L^1(M)$ to $L^1(N)$.
\par\vspace{3mm}
All elements of $L^1(M)$ can (and will) be interpreted as densely defined linear operators affiliated with $M$ acting on $L^2(M)$; full details are contained for instance in Appendix B of \cite{SS}. As $\Vert x^*\Vert_p=\Vert x\Vert_p$ for all $x\in M$, $p=1,2$, conjugation $x\mapsto x^*$ extends to an antilinear, isometric bijection of $L^p(M)$ onto itself. It is usually denoted by $J$ in the case $p=2$.
\par\vspace{3mm}
We denote by $\Tr$ the faithful, normal, semifinite trace on $\langle M,e_N\rangle$ characterized by 
$$
\Tr(xe_Ny)=\tau(xy)\quad\forall x,y\in M,
$$
and we will also make use of the so-called \emph{pull-down map} $\Phi$ defined first on the dense subalgebra $Me_NM:=\{\sum_{\mathrm{finite}} a_ie_Nb_i:\ a_i,b_i\in M\ \forall i\}$ by 
$$
\Phi(xe_Ny)=xy\quad\forall x,y\in M.
$$
It extends to an $M$-bimodule $*$-map from $L^1(\langle M,e_N\rangle,\Tr)$ to $L^1(M)$ and it has the following properties (see Section 4.5 in \cite{SS}, in particular Theorem 4.5.3):
\begin{enumerate}
	\item [(i)] $\tau(\Phi(x))=\Tr(x)$ for all $x\in L^1(\langle M,e_N\rangle,\Tr)$;
	\item [(ii)] $xe_N=\Phi(xe_N)e_N$ for all $x\in L^1(\langle M,e_N\rangle,\Tr)$;
	\item [(iii)] $\Vert\Phi(x)\Vert_1\leq\Vert x\Vert_{1,\Tr}$ for all $x\in L^1(\langle M,e_N\rangle,\Tr)$;
	\item [(iv)] $\Phi$ maps $Me_N$ onto $M$ with $\Vert\Phi(x)\Vert_2=\Vert x\Vert_{2,\Tr}$ for all $x\in Me_N$;
	\item [(v)] the equations $\Phi(xe_Ny)=xJy^*J\xi$ and $\Phi(ye_Nx)=Jx^*Jy\xi$ hold and both are in $L^2(M)$ for all all $x\in M$ and $y\in\langle M,e_N\rangle$.
\end{enumerate}
Moreover, we need the following extention of claim (iv) above: as $\Tr(e_N)=1$, one has $e_N\in L^2(\langle M,e_N\rangle,\Tr)$ as well, hence $Xe_N\in L^2(\langle M,e_N\rangle,\Tr)$ for every $X\in \langle M,e_N\rangle$, and $\Phi$ extends to $\langle M,e_N\rangle e_N$ isometrically, \emph{i.e.}
one has $Xe_N=\Phi(Xe_N)e_N$ with
$$
\Vert Xe_N\Vert_{2,\Tr}=\Vert \Phi(Xe_N)\Vert_2
$$
for every $X\in \langle M,e_N\rangle$. Indeed, it suffices to prove the above equality for $X=\sum_ia_ie_Nb_i$ where the sum is finite and $a_i,b_i\in M$ for every $i$. One has
$$
Xe_N=\sum_i a_ie_Nb_ie_N=\sum_i a_iE_N(b_i)e_N
$$
and 
$$
\Phi(Xe_N)=\sum_i a_iE_N(b_i)
$$
hence
\begin{eqnarray*}
\Vert Xe_N\Vert_{2,\Tr}^2 &=&
\Tr(e_NX^*Xe_N)\\
&=&
\Tr(\sum_{i,j}e_Nb_i^*e_Na_i^*a_je_Nb_je_N)\\
&=&
\sum_{i,j}\tau(E_N(b_i^*)E_N(a_i^*a_j)E_N(b_j))\\
&=&
\sum_{i,j}\tau(E_N(b_i^*)a_i^*a_jE_N(b_j))=\Vert \Phi(Xe_N)\Vert_2^2.
\end{eqnarray*}
 
\par\vspace{3mm}
Let $N$ be a $\Gamma$-invariant, unital von Neumann subalgebra of $M$. We denote by $u_\sigma(g)$ the unitary operator on $L^2(M)$ defined by $u_\sigma(g)x\xi=\sigma_g(x)\xi$ for every $x\in M$. Invariance implies that $E_N\sigma_g=\sigma_g E_N$, and that $u_\sigma(g)$ and $e_N$ commute for all $g\in \Gamma$. We extend $\sigma_g$ to an automorphism of $\langle M,e_N\rangle$ as follows: $\sigma_g(X)=u_\sigma(g) Xu_\sigma(g^{-1})$ for all $X\in \langle M,e_N\rangle$; it is characterized by $\sigma_g(ae_Nb)=\sigma_g(a)e_N\sigma_g(b)$ for all $a,b\in M$ and $g\in\Gamma$.
\par\vspace{3mm}
Moreover, there is a bijective correspondence between the set of projections of $\langle M,e_N\rangle$ and the right $N$-submodules of $L^2(M)$ given by $p\mapsto V=pL^2(M)$. The reciprocal map is denoted by $V\mapsto p_V$.
Note that $V$ is $u_\sigma(g)$-invariant if and only if $\sigma_g(p_V)=p_V$. Let $V$ be a right $N$-submodule of $L^2(M)$; following \cite{AET}, we say that $V$ is \emph{of finite lifted trace} if $\Tr(p_V)<\infty$. Furthermore, $V$ is said to be \emph{of finite rank} if there exist finitely many vectors $\xi_1,\ldots,\xi_r\in V$ such that 
$V=\overline{\sum_{i=1}^r \xi_iN}$. Following \cite{Popa}, we say that $\xi_1,\ldots,\xi_r$ is an \emph{orthonormal basis} of $V$ over $N$ if $E_N(\xi_i^*\xi_j)=\delta_{i,j}f_i$ where $f_i$ is a non zero projection for every $i=1,\ldots,r$. If it is the case, it follows immediately that, if $\eta\in V$ can be expressed as 
$$
\eta=\sum_{j=1}^r \xi_ja_j
$$
with $a_1=f_1a_1,\ldots,a_r=f_ra_r\in N$, then such a decomposition is unique: if $\eta=\sum_j\xi_jb_j$ with $b_j=f_jb_j\in N$ for all $j$, then $a_j=b_j$ for all $j=1,\ldots,r$.

\par\vspace{3mm}
We will also freely use comparison theory of projections in factors and its relationships with (semi)finite traces as it appears for instance in Part III of \cite{Dix} or in Chapter V of \cite{Tak}: we just recall that, if $\mathcal M$ is a von Neumann algebra, if $e,f\in\mathcal M$ are projections, then $e\preceq f$ if there exists a partial isometry $u\in\mathcal M$ such that $u^*u=e$ and $uu^*\leq f$. We denote by $\mathrm{Ip}(\mathcal M)$ the set of all partial isometries of $\mathcal M$.

\par\vspace{3mm}
Finally, if $\sigma:\Gamma\rightarrow \mathrm{Aut}(\mathcal M)$ is an action of $\Gamma$ on $\mathcal M$, a map $u:\Gamma\rightarrow \mathrm{Ip}(\mathcal M)$ is a $\sigma$-\emph{cocycle} if it satisfies $u(gh)=u(g)\sigma_g(u(h))$ for all $g,h\in\Gamma$.
\section{The case of subfactors}

Here is our first result.
Its proof is strongly inspired by that of Proposition 1.3 in \cite{PiPo}.

\begin{prop}
Suppose that $(M,\tau,\sigma,\Gamma)$ is a von Neumann dynamical system and that $N$ is a $\Gamma$-invariant type $\mathrm{II}_1$ subfactor of $M$. Let $p$ be a projection in $\langle M,e_N\rangle$ such that $0<t:=\Tr(p)<\infty$ and which is $\Gamma$-invariant. Then $V:=pL^2(M)$ has finite rank $r=t$ if $t$ is an integer, or
$r=n+1$
where $n$ is the integer part of $t$ otherwise, and there exists a family $(\xi_j)_{1\leq j\leq r}\subset L^2(M)$ such that:
\begin{enumerate}
	\item [(a)] $E_N(\xi_j^*\xi_k)=0$ for $j\not=k$;
	\item [(b)] $E_N(\xi_j^*\xi_j)=1$ for $1\leq j\leq n$, and $E_N(\xi_{n+1}^*\xi_{n+1}):=f_{n+1}$ is a projection of $N$ with trace equal to $\Tr(p)-n$;
	\item [(c)] each $\xi_je_N$ is a partial isometry in $\langle M,e_N\rangle$ and $\sum_j \xi_je_N\xi_j^*=p$;
	\item [(d)] every $\eta\in V$ has a unique decomposition $\eta=\sum_j \xi_jE_N(\xi_j^*\eta)$, thus in particular, the family 
	$(\xi_j)_{1\leq j\leq r}$ is an orthonormal basis of $V$ over $N$;
	\item [(e)] there exists a $\sigma\otimes i_r$-cocycle 
	$u: \Gamma  \rightarrow  \mathrm{Ip}(M_r(N))$,
	$g  \mapsto  u(g)=(u_{i,j}(g))_{1\leq i,j\leq r}$, 
	such that
	$$
	u_\sigma(g)\xi_j=\sum_i \xi_iu_{i,j}(g)\quad\forall j,g.
	$$
\end{enumerate}
\end{prop}
\textsc{Proof.} Let us assume that $\Tr(p)$ is not an integer, hence that $0<\Tr(p)-n<1$. Since $\langle M,e_N\rangle$ is a factor of type II, there exist projections $g_1,\ldots,g_n,g_{n+1}\in\langle M,e_N\rangle$ such that 
\begin{enumerate}
	\item [(i)] $g_jg_k=0$ if $j\not=k$;
	\item [(ii)] $\sum_j g_j=p$;
	\item [(iii)] $\Tr(g_j)=\Tr(e_N)=1$ for all $j\leq n$, and $\Tr(g_{n+1})=\Tr(p)-n$.
\end{enumerate}
Hence, for every $j\leq n$, $g_j$ is a projection in $\langle M,e_N\rangle$ which is equivalent to $e_N$ (since they have the same trace, and since $JN'J=\langle M,e_N\rangle$ is a factor), and $g_{n+1}\preceq e_N$. Thus, there exist partial isometries $v_j\in \langle M,e_N\rangle$, $j=1,\ldots,n+1$, such that $v_jv_j^*=g_j$ for all $j=1,\ldots,n+1$, $v_j^*v_j=e_N$ for all $j\leq n$, and $v_{n+1}^*v_{n+1}\leq e_N$.
\par
Since $v_j=v_je_N$ for all $j$, the use of the pull-down map $\Phi$ implies that $v_j=\xi_je_N$ with $\xi_j=\Phi(v_j)=\Phi(v_je_N)\in L^2(M)$, and this proves immediately statement (c).
\par
For $j\not=k$, one has
$$
0=v_j^*v_k=e_Nv_j^*v_ke_N=E_N(\xi_j^*\xi_k)e_N
$$
hence $E_N(\xi_j^*\xi_k)=0$, which proves statement (a). Similarly, the fact that $v_j^*v_j=e_N$ for $j\leq n$ 
shows that $E_N(\xi_j^*\xi_j)=1$ for $j\leq n$. 
If $j=n+1$, then $v_{n+1}^*v_{n+1}=E_N(\xi_{n+1}^*\xi_{n+1})e_N$ is a subprojection of $e_N$, hence $E_N(\xi_{n+1}^*\xi_{n+1}):=f_{n+1}$ is a projection, too. Thus claim (b) holds true.
\par
If $x\in M$, one has
$$
p(x\xi)=p(xe_N\xi)=\sum_j\xi_je_N\xi_j^*xe_N\xi=\sum_j\xi_jE_N(\xi_j^*x\xi)
$$
thus, by density of $M$ in $L^2(M)$, this proves statement (d).
\par
Finally, put $u_{i,j}(g)=E_N(\xi_i^*u_\sigma(g)\xi_j)\in L^1(N)$ for all $i,j$, so that 
$$
u_\sigma(g)\xi_j=\sum_i \xi_i u_{i,j}(g)\quad\forall j,g
$$
and that $u_{i,j}(g)e_N=e_Nu_{i,j}(g)e_N=e_NE_N(\xi^*u_\sigma(g)\xi_j)e_N$ for all $i,j$ and $g$.
Then we claim that $\sum_k u_{k,i}(g)^*u_{k,j}(g)=\delta_{i,j}\sigma_g E_N(\xi_i^*\xi_j)$; this will prove that all $u_{i,j}(g)$'s are bounded operators (by taking $i=j$), and that $u(g)\in \mathrm{Ip}(M_r(M))$. One has, indeed,
\begin{eqnarray*}
\sum_k u_{k,i}(g)^*u_{k,j}(g)e_N 
&=&
\sum_k e_Nu_{k,i}(g)^*u_{k,j}(g)e_N\\
&=&
\sum_k e_N(u_\sigma(g)\xi_i)^*\xi_k e_N\xi_k^* u_\sigma(g)\xi_j e_N\\
&=&
e_N(u_\sigma(g) \xi_i)^* pu_\sigma(g) \xi_j e_N\\
&=&
\delta_{i,j}\sigma_g(E_N(\xi_i^*\xi_j))e_N
\end{eqnarray*}
since $u_\sigma(g)$ commutes with $p$ and $p\xi_i=\xi_i$. This shows that $u_{k,i}(g)$ is bounded, hence belongs to $N$ and that 
$$
u(g)^*u(g)=\left(
\begin{array}{cccc}
1 & 0 & \ldots & 0\\
0 & 1 & \ldots & 0\\
\vdots & \ddots & \ddots & \vdots\\
0 & 0 & \ldots & \sigma_g(f_{r})
\end{array}\right),
$$
proving that $u(g)$ is a partial isometry for all $g$. As obviously $f_{r}u_{r,r}(g)\sigma_g(f_r)=u_{r,r}(g)$ for all $g$, the decomposition of each $u_\sigma(g)\xi_j$ is unique, we see that $u(gh)=u(g)\sigma_g\otimes i_r(u(h))$ for all $g,h\in\Gamma$ as follows:
\begin{eqnarray*}
u_\sigma(gh)\xi_j 
&=&
u_\sigma(g)(u_\sigma(h)\xi_j)=u_\sigma(g)\left(\sum_i \xi_i u_{i,j}(h)\right)\\
&=&
\sum_i u_\sigma(g)(\xi_i u_{i,j}(h))=
\sum_i u_\sigma(g)\xi_i\sigma_g(u_{i,j}(h))\\
&=&
\sum_k \xi_k \left(\sum_i \xi_i u_{k,i}(g)\sigma_g(u_{i,j}(h))\right)\\
&=&
\sum_k \xi_k u_{k,j}(gh).
\end{eqnarray*}
\hfill $\square$
\par\vspace{3mm}\noindent
\textbf{Remark.} Contrary to what is stated in Lemma 4.1 of \cite{AET}, the \emph{entries} $u_{i,j}$ of the unitary operator $u$ are \emph{not} unitary operators themselves (neither are the entries $u_{i,j}^{(-n)}$ in the proof of Proposition 4.5), but they satisfy $\Vert u_{i,j}\Vert\leq 1$ for all $i,j$, and this suffices to get the conclusion in the proof of Proposition 4.5 in \cite{AET}.

\section{The general case}

Here is our last result.

\begin{prop}
Let $(M,\tau,\sigma,\Gamma)$ be a von Neumann dynamical system, let $N$ be a unital, $\Gamma$-invariant von Neumann subalgebra of $M$, let $p$ be a projection in $\langle M,e_N\rangle$ such that $0<\Tr(p)<\infty$ and $\sigma_g(p)=p$ for all $g\in\Gamma$. For every $\varepsilon>0$, there exists a projection $q\leq p$ in $\langle M,e_N\rangle$ such that $\sigma_g(q)=q$ for all $g$, $\Tr(p-q)<\varepsilon$
and there exists a finite family $(\xi_j)_{1\leq j\leq r}\subset L^2(M)$ with the following properties:
\begin{enumerate}
	\item [(a)] $E_N(\xi_j^*\xi_k)=0$ for all $j\not=k$;
	\item [(b)] $E_N(\xi_j^*\xi_j)=:f_j$ is a non zero projection of $N$ for all $j=1,\ldots,r$;
	\item [(c)] each $\xi_je_N$ is a partial isometry and $\sum_j \xi_je_N\xi_j^*=q$;
	\item [(d)] every $\eta\in qL^2(M)$ has a decomposition $\eta=\sum_j \xi_jE_N(\xi_j^*\eta)$, thus in particular, 
	$W:=qL^2(M)$ is a right $N$-submodule of finite rank and
	the family $(\xi_j)_{1\leq j\leq r}$ is an orthonormal basis of $W$ over $N$;
	\item [(e)] there exists a $\sigma\otimes i_r$-cocycle 
	$u: \Gamma  \rightarrow  \mathrm{Ip}(M_r(N))$,
	$g  \mapsto  u(g)=(u_{i,j}(g))_{1\leq i,j\leq r}$, 
	such that
	$$
	u_\sigma(g)\xi_j=\sum_i \xi_iu_{i,j}(g)\quad\forall j,g.
	$$
\end{enumerate}
\end{prop}
\textsc{Proof.} Assume first that $p\preceq e_N$. Then there exists a partial isometry $v\in\langle M,e_N\rangle$ such that $vv^*=p$ and $v^*v\leq e_N$. Thus $v=ve_N=\xi_1e_N$ where $\xi_1=\Phi(v)=\Phi(ve_N)\in L^2(M)$. One checks as in the proof of Proposition 1 that $v^*v=E_N(\xi_1^*\xi_1)=f$ is a projection in $N$
and that
$p\eta=\xi_1E_N(\xi_1^*\eta)$ for every $\eta\in L^2(M)$, hence $pL^2(M)=\overline{\xi_1N}$. 
\par
Moreover, if we set $u(g)=E_N(\xi_1^*u_\sigma(g)\xi_1)$, then it is easy to see that $u(g)^*u(g)e_N=\sigma_g(f)e_N$, hence that $u(g)\in N$ is a partial isometry and one can check as in the proof of the preceeding proposition that $u$ is a cocycle.
\par
Next, if $p\not\preceq e_N$, by the Comparability Theorem (see Theorem V.1.8 of \cite{Tak}, for instance), there exists a central projection $z\in\langle M,e_N\rangle$ such that $e_Nz\preceq pz$ and $p(1-z)\preceq e_N(1-z)$. By our present assumption, $z\not=0$ and $e_Nz\not=0$, since the central support of $e_N$ is 1. Then we define the following set $\mathcal Z$, ordered by inclusion: it is the set of all families of projections $(g_i)_{i\in I}\subset \langle M,e_N\rangle$ such that $g_ig_j=0$ for all $i\not=j$, $0\not=g_i\preceq e_N$, $\sigma_g(g_i)=g_i$ for all $i\in I$ and all $g\in\Gamma$ and finally $\sum_{i\in I}g_i\leq p$. 
\par
Let us see that $\mathcal Z\not=\emptyset$: as $e_Nz\preceq pz$, there exists a partial isometry $u\not=0$ such that
$$
u^*u= e_Nz\leq e_N\quad\textrm{and}\quad
uu^*\leq pz\leq p.
$$
Then $e:=\bigvee_{g\in\Gamma}\sigma_g(uu^*)$ is a $\Gamma$-invariant non zero projection, $e\leq p$ and, by Lemma V.2.2 of \cite{Tak}, $e\preceq e_N$. Hence, by Zorn's Lemma, let $(g_i)_{i\in I}$ be a maximal element of $\mathcal Z$.
\par
 We claim that 
$\sum_{i\in I}g_i=p$. Indeed, if it is not the case, put $f=p-\sum_ig_i\not=0$. Then $\sigma_g(f)=f$ for all $g$, $f\leq p$, and, as the central support of $e_N$ is 1, there exists a non zero projection $f'\leq f$ such that $f'\preceq e_N$ (Lemma V.1.7 of \cite{Tak}). As above, the projection $f''=\bigvee_{g\in \Gamma}\sigma_g(f')$ is non zero, $\Gamma$-invariant, and the family $(g_i)_{i\in I}\cup\{f''\}$ contradicts the maximality of $(g_i)_{i\in I}$. 
\par
As $\sum_i \Tr(g_i)=\Tr(p)<\infty$, if $\varepsilon>0$ is given, there exists $\{1,\ldots,r\}\subset I$ such that
$\Tr(p-q)<\varepsilon$ and $\sigma_g(q)=q$ for all $g$, where we have set $q=\sum_{i=1}^r g_i$. Now, as in the proof of Proposition 1, there exist partial isometries $v_1,\ldots,v_r\in\langle M,e_N\rangle$ such that $v_iv_i^*=g_i$ and $v_i^*v_i\leq e_N$ for all $i$. One can find vectors $\xi_1,\ldots,\xi_r\in L^2(M)$ such that $v_i=\xi_ie_N$ for all $i$, and it is easy to see that statements (a), (b), (c) and (d) hold true. Furthermore, if we set $u_{i,j}(g)=E_N(\xi_i^*u_\sigma(g) \xi_j)$ so that $u_\sigma(g)\xi_j=\sum_i\xi_i u_{i,j}(g)$, 
one can prove as in the previous section that 
$$
\sum_k u_{k,i}(g)^*u_{k,j}(g)e_N=\delta_{i,j}\sigma_g(f_i)e_N
$$
for all $i,j$, proving that $u(g)=(u_{i,j}(g))\in M_r(N)$ is a partial isometry. The cocycle identity is proved exactly as in Proposition 1.
\hfill $\square$

\section{An application to relative weak mixing}

Let us recall Definition 3.7 of \cite{AET}:

\begin{defn}
Let $(M,\tau,\alpha)$ be a von Neumann dynamical system and let $N\subset M$ be an $\alpha$-invariant von Neumann subalgebra. Then $\alpha$ is \textbf{weakly mixing relative to} $N$ if, for any $a\in M\ominus N$ we have
\[
\lim_{N\to\infty}\frac{1}{N}\sum_{n=1}^N \Vert E_N(a^*\alpha^n(a))\Vert^2_2=0.
\]
\end{defn}

Before stating the main result of the present section, we need to recall well known equivalent conditions on bounded sequences of complex numbers (see for instance \cite{Walt}, Theorem 1.20): First, recall that a set $Z\subset\N$ has \textbf{zero density} if 
\[
\lim_{n\to\infty}\frac{|Z\cap\{0,1,\ldots,n-1\}|}{n}=0.
\]

\begin{prop}
(\cite{Walt}, Theorem 1.20)
For a bounded sequence $(a_n)_{n\geq 0}\subset \C$, the following conditions are equivalent:
\begin{enumerate}
\item [(i)] one has $\displaystyle{\lim_{N\to\infty}\frac{1}{N}\sum_{n=0}^{N-1}|a_n|=0}$;
\item [(ii)] there exists $Z\subset\N$ with zero density such that
$\displaystyle{
\lim_{n\to\infty,n\notin Z} a_n=0}$;
\item [(iii)] one has 
$\displaystyle{\lim_{N\to\infty}\frac{1}{N}\sum_{n=0}^{N-1}|a_n|^2=0}$.
\end{enumerate}
\end{prop}

The previous proposition will be used in the proof of our last result.

\begin{prop}
Let $(M,\tau,\alpha)$ be a von Neumann dynamical system and let $N\subset M$ be an $\alpha$-invariant von Neumann subalgebra of $M$. Then the following two conditions are equivalent:
\begin{enumerate}
\item [(1)] For every finite set $F\subset L^2(M)\ominus L^2(N)$ such that $E_N(\eta^*\eta)$ is bounded for every $\eta\in F$ and for every $\ep>0$, there exists $n>0$ such that 
\[
\max_{\eta,\zeta\in F}\Vert E_N(\eta^*u_\alpha^n\zeta)\Vert_2\leq \ep.
\]
\item [(2)] The automorphism $\alpha$ is weakly mixing relative to $N$. 
\end{enumerate}
\end{prop}
\textsc{Proof.} $(1)\ \Rightarrow\ (2)$: If (2) does not hold true, Proposition 3.8 in \cite{AET} implies the existence of a $u_\alpha$-invariant right submodule $V\subset L^2(M)\ominus L^2(N)$ such that $\Tr(p_V)<\infty$. Thus, Proposition 2 implies that, if $\ep>0$ is small enough, there exists a $u_\alpha$-invariant right submodule $W\subset V$ and an orthonormal basis $\xi_1,\ldots,\xi_r\in W$ such that $E_N(\xi_i^*\xi_j)=\delta_{i,j}f_j$ for all $i,j$, where, $f_j\in N$ is a non-zero projection for every $j$. Moreover, there exists a cocycle $\omega: \{\alpha^k\colon k\in \Z\}\rightarrow M_r(N)$ such that
\[
u_\alpha^n\xi_j=\sum_{i=1}^r \xi_i\omega_{i,j}(\alpha^n)\quad \forall n,\forall 1\leq j\leq r.
\]
One also has
\[
\sum_{k=1}^r \omega_{k,i}(\alpha^n)^*\omega_{k,j}(\alpha^n)=\delta_{i,j}\alpha^n(f_i) \quad \forall i,j,n.
\]
As $E_N(\xi_i^*u_\alpha^n\xi_j)=\omega_{i,j}(\alpha^n)$, one gets
\begin{align*}
\sum_{i,j=1}^r \Vert E_N(\xi_i^*u_\alpha^n\xi_j)\Vert_2^2
&=
\sum_{i,j=1}^r \tau(\omega_{i,j}(\alpha^n)^*\omega_{i,j}(\alpha^n))\\
&=
\sum_{j=1}^r\tau\Big(\sum_{i=1}^r \omega_{i,j}(\alpha^n)^*\omega_{i,j}(\alpha^n)\Big)\\
&=
\sum_{j=1}^r \tau(\alpha^n(f_j))=\sum_{j=1}^n \tau(f_j)\eqqcolon \delta>0
\end{align*}
for every $n\in \Z$. \\
If $\alpha$ satisfied (1), taking $F=\{\xi_1,\ldots,\xi_r\}$ and $\displaystyle{\ep=\frac{\sqrt \delta}{\sqrt 2 r}}$, there would exists an $n$ such that 
\[
\max_{i,j}\Vert E_N(\xi_i^*u_\alpha^n\xi_j)\Vert_2^2\leq \frac{\delta}{2r^2}
\]
hence 
\[
\delta\leq \sum_{i,j=1}^r \Vert E_N(\xi_i^*u_\alpha^n\xi_j)\Vert_2^2\leq \delta/2\]
which yields a contraction.\\
$(2)\ \Rightarrow\ (1)$: By Proposition 4, condition (2) implies that, for every finite set $F\subset M\ominus N$, there exists $Z\subset \N$ with zero density such that
\begin{equation*} 
\lim_{n\to\infty, n\notin Z}\Vert E_N(a^*\alpha^n(a))\Vert_2=0
\end{equation*}
for every $a\in F$. By density, it suffices to verify (2) for $F\subset M\ominus N$ finite. But the following polar decomposition holds true:
\[
E_N(x^*\alpha^n(y))=\frac{1}{4}\sum_{k=0}^3 i^k E_N((y+i^kx)^*\alpha^n(y+i^kx)).
\]
Applying the above limit to all elements of $F'\coloneqq \{y+i^kx\colon x,y\in F, 0\leq k\leq 3\}\subset M\ominus N$, we get (1).
\hfill $\square$

\par\vspace{1cm}
\bibliographystyle{plain}
\bibliography{max}
\par
\vspace{1cm}
\noindent
\begin{flushright}
     \begin{tabular}{l}
       Universit\'e de Neuch\^atel,\\
       Institut de Math\'emathiques,\\       
       Emile-Argand 11\\
       CH-2000 Neuch\^atel, Switzerland\\
       \small {pajolissaint@gmail.com}
     \end{tabular}
\end{flushright}

\end{document}